\documentclass[iccmp]{ipbook}
\usepackage{amsmath, amsthm, amssymb, amsfonts}
\usepackage{mathrsfs}
\usepackage{cleveref}
\usepackage{graphicx}
\usepackage{subcaption}
\usepackage{diagbox}
\usepackage{xcolor}
\usepackage{empheq}

\startlocaldefs
\endlocaldefs

\numberwithin{equation}{section}

\theoremstyle{definition}
\newtheorem{theorem}{Theorem}[section]
\newtheorem{definition}[theorem]{Definition}

\DeclareMathOperator*{\argmax}{arg\,max}
\DeclareMathOperator*{\argmin}{arg\,min}

\usepackage{tikz}
\usetikzlibrary{arrows.meta, positioning}
\usetikzlibrary{shapes.geometric, arrows, positioning}
\usepackage{pgfplots}
\usetikzlibrary{calc,intersections}

\tikzstyle{arrow} = [->, >=stealth, -triangle 60]

\usetikzlibrary{shapes.geometric, arrows, positioning}

\tikzset{
	mybox/.style  = {draw, rectangle, minimum width=4cm, minimum height=0.8cm, text centered, text width=4.4cm,   
		font=\normalsize},
	box/.style  = {draw, rectangle, minimum width=2.0cm, minimum height=0.6cm, text centered, text width=3.0cm,   
		font=\normalsize},
	myarrow/.style = {line width=0.2pt, draw=black, -triangle 60, postaction={draw, line width=0.2pt, shorten >=10pt,-}}
}

\firstpage{1}
\lastpage{1}

\title[Modern Theory of Gradient-Based Optimization]{Modern Theory of Gradient-Based Optimization}
\author[Bin Shi]{Bin Shi}
\address{Center for Mathematics and Interdisciplinary Sciences, Fudan University, Shanghai 200433, China}
\address{Shanghai Institute for Mathematics and Interdisciplinary Sciences, Shanghai 200433, China}
\email{binshi@fudan.edu.cn}\thanks{This work was partially supported by the NSFC (Grant No. 12241105) and by SIMIS (startup fund and cross-disciplinary research projects).}

\begin{document}

\begin{abstract}
In this review, we offer a comprehensive survey of emerging techniques in gradient-based optimization, with a particular emphasis on the interplay between ordinary differential equation (ODE) perspectives and their extensions into discrete Lyapunov analysis. We begin by examining the acceleration mechanisms underlying Nesterov’s accelerated gradient method for strongly convex functions (\texttt{NAG-SC}) and Polyak’s heavy-ball method, identifying the gradient-correction term as the primary driver of acceleration. This mechanistic insight is substantiated through high-resolution ODE modeling and the systematic construction of Lyapunov functions. We then synthesize recent advancements in convex optimization regarding \texttt{NAG} and its proximal generalization, the fast iterative shrinkage-thresholding algorithm (\texttt{FISTA}). Key topics include the accelerated convergence of gradient norms, underdamped acceleration, linear convergence under strong convexity, and novel Lyapunov frameworks for establishing convergence and monotonicity properties of generalized accelerated methods. Furthermore, we demonstrate how these ODE approximations and Lyapunov techniques can be extended to provide a unified framework for analyzing advanced optimization algorithms, including the alternating direction method of multipliers (\texttt{ADMM}), the primal-dual hybrid gradient (\texttt{PDHG}) method, and their respective accelerated variants. Finally, we discuss recent progress in minimax optimization and outline future directions for extending Lyapunov-based analysis to saddle-point problems.
\end{abstract}

%First, the acceleration mechanisms of Nesterov’s accelerated gradient method (\texttt{NAG-SC}) and Polyak’s heavy-ball method are examined, where the gradient-correction term is identified as the primary driver of acceleration. This perspective is supported by high-resolution ODEs and systematic constructions of Lyapunov functions. Next, recent advancements in \texttt{NAG} and its proximal generalization, the fast iterative shrinkage-thresholding algorithm (\texttt{FISTA}), are reviewed; specifically, topics such as accelerated convergence of gradient norms, underdamped acceleration, linear convergence under strong convexity, and new Lyapunov frameworks for monotonic variants are addressed. Subsequently, it is demonstrated how these ODE approximations and Lyapunov techniques provide a unified framework for analyzing advanced methods, including the Alternating Direction Method of Multipliers (\texttt{ADMM}), the Primal-Dual Hybrid Gradient (\texttt{PDHG}) method, and their accelerated variants. Finally, recent progress on minimax optimization is discussed, and future directions for extending Lyapunov-based analysis to saddle-point problems are outlined.

\maketitle

%\tableofcontents

%%  The body
\section{Introduction}
\label{sec: intro}

The study of gradient-based optimization algorithms has a storied history, perhaps tracing its origins back to the pioneering work of Cauchy~\cite{cauchy1847methode}. However, the quest for truly ``accelerated'' methods gained significant momentum only in the 1960s, largely inspired by the insights of Gelfand and Tsetlin~\cite{gelfand1961printsip}. Building upon these advancements, Polyak~\cite{polyak1964some} bridged the gap between optimization and classical mechanics. Drawing from invariant manifold theory, a local theory concerning equilibria of dynamical system, he proposed the heavy-ball method. This contribution effectively introduced the concept of inertia into the optimization process, allowing iterates to ``roll'' past local oscillations. The field underwent a paradigm shift when Nemirovski and Yudin~\cite{nemirovsky1983problem} established computational complexity as a formal framework for evaluating algorithmic performance.  This theoretical rigor set the stage for Nesterov, who developed the ``Estimate Sequence'' technique,  which led to the landmark Nesterov's accelerated gradient (\texttt{NAG}) method for convex objective functions~\cite{nesterov1983method}. He later refined this approach for strongly convex functions~(\texttt{NAG-SC}), with a comprehensive mathematical treatment presented in~\cite{nesterov2018lectures}.

Since the turn of this century, machine learning has experienced explosive growth, largely propelled by the development of efficient optimization algorithms. A major bottleneck in large-scale optimization is the substantial computational and memory cost required to evaluate or store Hessian matrices, which often renders second-order methods impractical for modern high-dimensional datasets. Consequently, gradient-based methods, valued for their low per-iteration overhead, have become the engine driving recent advancements. With the rise of deep learning, \texttt{NAG} has become widely regarded as a momentum-based method of central importance~\cite{sutskever2013importance}. To better understand its behavior, a perspective that investigates~\texttt{NAG} through the lens of ordinary differential equations (ODEs) is introduced in~\cite{su2016differential}, using a continuous-time approximation to model its dynamics. Departing from earlier local approaches~\cite{polyak1964some}, this framework adopts a numerical analysis standpoint centered on discretization error. When combined with the Lyapunov function technique for global analysis, as proposed in~\cite{su2016differential}, this method opens up a new methodology for studying optimization algorithms. While the technical seeds of this framework can be traced back to~\cite{su2016differential}, its true utility as a robust and essential analytical tool was first recognized and fully articulated in Jordan's ICM 2018 plenary lecture~\cite{jordan2018dynamical}.

This review begins by adopting the perspective of high-resolution ODEs and the discrete Lyapunov analysis introduced in~\cite{shi2022understanding} to address the fundamental question proposed in~\cite{jordan2018dynamical}:``\textbf{What is the underlying mechanism that generates acceleration?}'' Building upon this foundation,~\Cref{sec: acceleration} provides  a detailed analysis for the acceleration mechanism.  The framework’s utility is then demonstrated through broader applications:~\Cref{sec: nag-fista} introduces several new results for~\texttt{NAG} and its proximal version, the Faster Iterative Shrinkage-Thresholding Algorithm (\texttt{FISTA}). Furthermore,~\Cref{sec: admm-pdhg} extends this analytical lens to other prominent algorithms, such as the Alternating Direction Method of Multipliers (\texttt{ADMM}) and the Primal-Dual Hybrid Gradient (\texttt{PDHG}) methods, including their various extensions.  Finally,~\Cref{sec: conclusion} provides concluding remarks and outlines potential directions for future research.

%The differential equation as a continuous approximation as a perspective to investigate the continuous dynamics has introduced to investigate the, which is introduced in. Different from the previous work, here introduces the numerical scheme  and the Lyapunov function technique has been introduced for the global analysis instead of the local theory. The real power of this technique has been recognized and as an important technique proposed in

\section{The acceleration mechanism}
\label{sec: acceleration}

In this section, we assume that the objective function $f$ is sufficiently smooth and  satisfies both $L$-smoothness and $\mu$-strong convexity. 

\begin{definition}
\label{defn: strongly-convex}
For any $x, y \in \mathbb{R}^n$, there exists two constants $L, \mu > 0$ such that the following conditions hold: 
\begin{itemize}
\item \textbf{$L$-smoothness}: The gradient $\nabla f$ is Lipschitz continuous with constant $L$:
\begin{equation}
\label{eqn: l-smooth}
\left\| \nabla f(x) - \nabla f(y) \right\| \leq L \| x - y \|. 
\end{equation}

\item \textbf{$\mu$-strong convexity}: The function $f$ satisfies: 
\begin{equation}
\label{eqn: mu-strongly}
f(y) \geq f(x) + \left\langle  \nabla f(x), y - x \right\rangle + \frac{\mu}{2} \|  y - x \|^2.
\end{equation}
\end{itemize}
The class of objective functions satisfies~\eqref{eqn: l-smooth} and~\eqref{eqn: mu-strongly} is denoted as $\mathcal{S}^{1}_{\mu, L}$. 
\end{definition}
According  \Cref{defn: strongly-convex}, any function $f \in \mathcal{S}^{1}_{\mu, L}$ can be equivalently expressed by the following inequalities for all $x, y \in \mathbb{R}^n$: 
\begin{itemize}
\item \textbf{Gradient variation bounds}:
\begin{equation}
\label{eqn: l-smooth-mu-strong-equiv-1}
\mu  \| x - y \| \leq \left\| \nabla f(x) - \nabla f(y) \right\| \leq L \| x - y \|. 
\end{equation}
\item \textbf{Quadratic growth bounds}: 
\begin{equation}
\label{eqn: l-smooth-mu-strong-equiv-2}
\frac{\mu}{2}  \| x - y \|^2 \leq f(y) - f(x) - \left\langle  \nabla f(x), y - x \right\rangle  \leq \frac{L}{2} \| x - y \|^2. 
\end{equation}
\end{itemize}

With the properties~\eqref{eqn: l-smooth} ---~\eqref{eqn: l-smooth-mu-strong-equiv-2} established,  we consider standard gradient-based optimization algorithms. Given a step size $0 < s \leq 1/L$ and an initial point $x_0 \in \mathbb{R}^d$,  the gradient descent update is given by 
\begin{equation}
\label{eqn: gd}
x_{k+1} = x_{k} - s \nabla f(x_k).
\end{equation}
Taking the gradient descent method~\eqref{eqn: gd} as a baseline,  we introduce the heavy-ball method. For any initial $x_0, x_1 \in \mathbb{R}^d$, the update rule incorporates a momentum term as:
\begin{equation}
\label{eqn: hb}
x_{k+1} = x_{k} - s \nabla f(x_k) + \frac{1 - \sqrt{\mu s}}{1 + \sqrt{\mu s}} (x_{k} - x_{k-1}).
\end{equation}
Finally,~\texttt{NAG-SC} is defined by the following coupled iterations, starting with $x_0 = y_0 \in \mathbb{R}^d$:
\begin{equation}
\label{eqn: nag-sc}
\left\{ \begin{aligned}
         & y_{k+1} = x_{k}  - s \nabla f(x_k), \\
         & x_{k+1} = y_{k+1} + \frac{1 - \sqrt{\mu s}}{1 + \sqrt{\mu s}} \left( y_{k+1} - y_{k} \right).
         \end{aligned} \right. 
\end{equation}
If we expand~\texttt{NAG-SC} solely in terms of the iterate $\{ x_{k} \}_{k=0}^{\infty}$, we obtain: 
\begin{align}
x_{k+1} = x_{k} - s \nabla f(x_k) & + \frac{1 - \sqrt{\mu s}}{1 + \sqrt{\mu s}} (x_{k} - x_{k-1}) \nonumber \\ 
                                                   & - \underbrace{\frac{1 - \sqrt{\mu s}}{1 + \sqrt{\mu s}} \cdot s \left( \nabla f(x_{k}) - \nabla f(x_{k-1}) \right)}_{\mathbf{Gradient \; Correction}}, \label{eqn: nag-sc-x} 
\end{align}
which is known as the gradient-correction scheme, discussed extensively in \cite{shi2022understanding}. Alternatively, expanding~\texttt{NAG-SC} in terms of the iterate $\{ y_{k} \}_{k=0}^{\infty}$ yields: 
\begin{equation}
\label{eqn: nag-sc-y}
y_{k+1} = y_{k} - s \nabla f\bigg(  y_{k} +  \underbrace{\frac{1 - \sqrt{\mu s}}{1 + \sqrt{\mu s}} \left( y_{k} - y_{k-1} \right)}_{\mathbf{Implicit\; Velocity}} \bigg)  + \frac{1 - \sqrt{\mu s}}{1 + \sqrt{\mu s}} (y_{k} - y_{k-1}).
\end{equation}
where the difference $y_{k} - y_{k-1}$ can be viewed as the velocity, so we refer it to as the implicit-velocity scheme, as discussed in \cite{chen2025revisiting}.

\begin{figure}[htp!]
\centering
\includegraphics[scale=0.22]{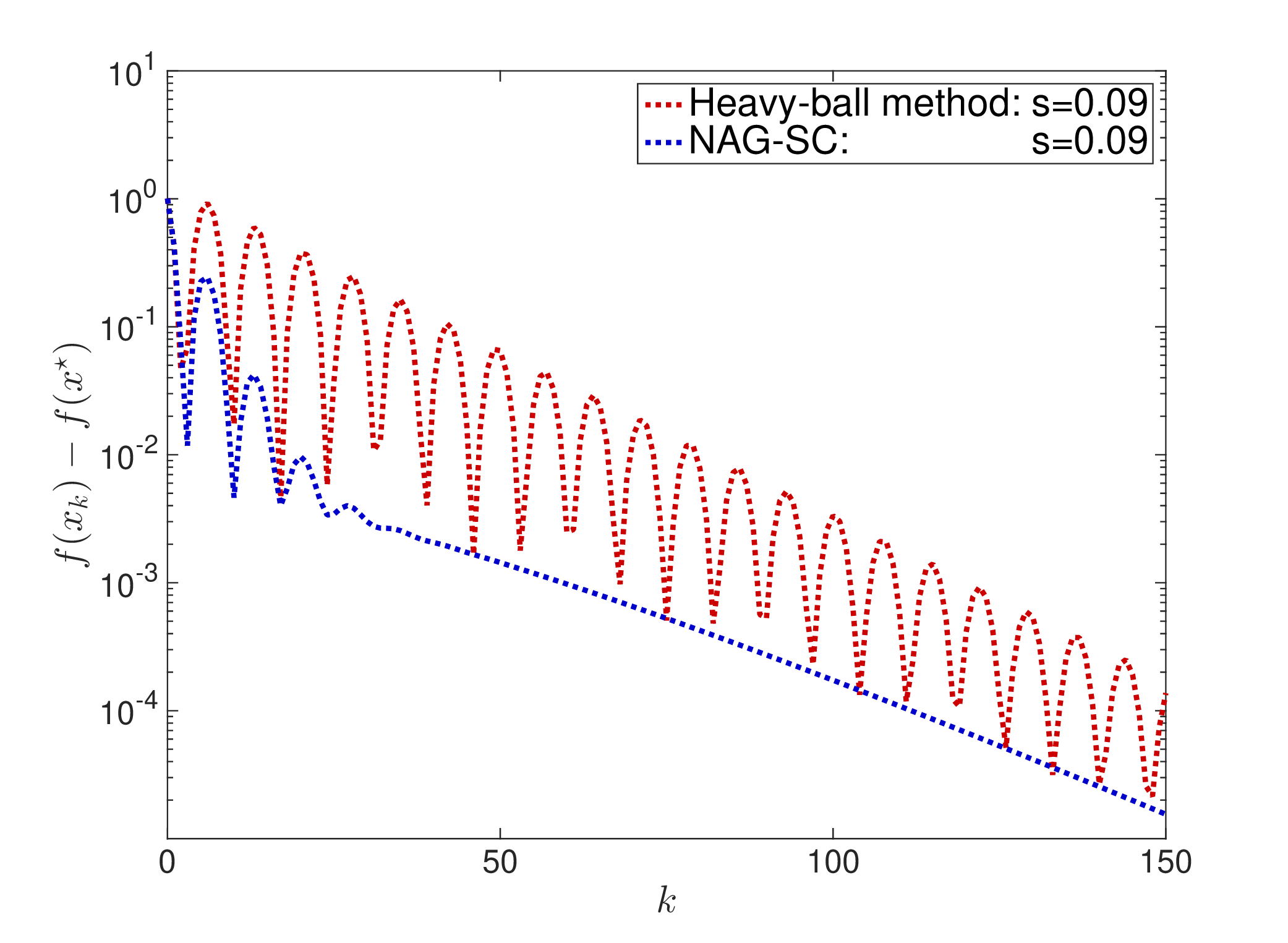}
\caption{Comparison of the convergence behaviors of the heavy-ball method and~\texttt{NAG-SC}.} 
\label{fig: nag-sc-and-heavy}
\end{figure}

Although these three gradient-based algorithms share closely related iterative structure, their convergence rates differ significantly, as summarized in~\Cref{tab: convergence-rate}. 
\begin{table}[htp!]
\centering
\scalebox{0.8}{
 \begin{tabular}{|l|l|} 
   \hline
\diagbox{Algorithms}{Objective}& $\mu$-Strongly Convex Objective  \\ 
  %\midrule 
  \hline
  Gradient Descent~\eqref{eqn: gd}                  & \textbf{Global}                            $\text{O}\left( \left(1 - \frac{\mu}{L} \right)^k \right)$                   \\
  \hline
  Heavy-Ball Method~\eqref{eqn: hb}                & \textbf{Local}                              $\text{O}\left( \left(1 - \sqrt{\frac{\mu}{L} }\right)^k \right)$         \\
  \hline
  \texttt{NAG-SC}~\eqref{eqn: nag-sc}              &  \textbf{Global}                           $\text{O}\left( \left(1 - \sqrt{\frac{\mu}{L} }\right)^k \right)$          \\ 
  \hline
 \end{tabular} 
 }
\caption{Convergence rates of gradient-based algorithms.} 
\label{tab: convergence-rate}
\end{table}
This discrepancy raises a fundamental question: \textbf{Why does~\texttt{NAG-SC}~\eqref{eqn: nag-sc} achieve global acceleration, whereas the heavy-ball method~\eqref{eqn: hb} only accelerates locally?} A direct comparison suggests that this distinction arises from the effects of ``gradient-correction'' or ``implicit-velocity''. To further investigate the mechanism, we first illustrate their convergence behaviors in~\Cref{fig: nag-sc-and-heavy}.  We observe that, while the heavy-ball method~\eqref{eqn: hb} oscillates throughout the convergence process,~\texttt{NAG-SC}~\eqref{eqn: nag-sc} exhibits only minor initial oscillations before transitioning into a smooth and stable convergence regime.

To understand the acceleration phenomenon, a high-resolution ODE framework was introduced in \cite{shi2022understanding} and subsequently refined in \cite{chen2025revisiting}. As  visually outlined in~\Cref{fig:chart}, this framework consists of the following four steps: 

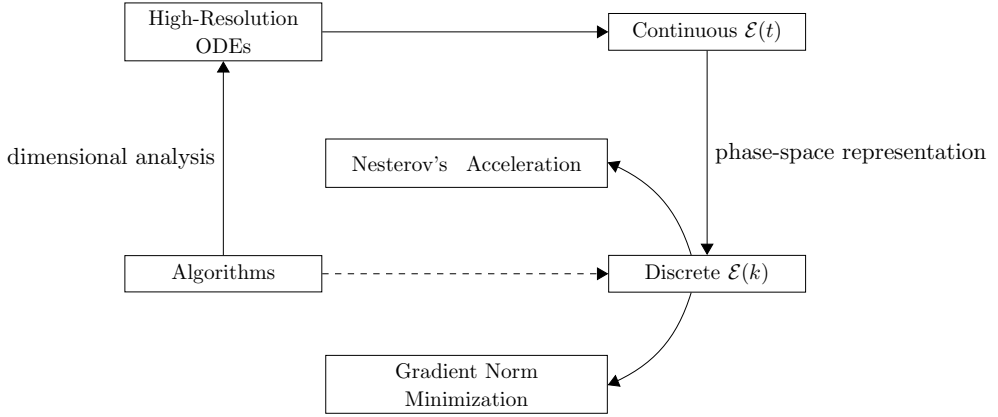
\begin{figure}[!htp]
\begin{center}
\scalebox{0.8}{
\begin{tikzpicture}[node distance=4cm]
    \node (n00) [box, draw=black] {\small Algorithms};
    \node (n10) [box, draw=black,  above of=n00] {\small High-Resolution ODEs};
    \node (n11) [box, draw=black,  above of=n00, xshift=+8cm] {\small Continuous $\mathcal{E}(t)$};
    \node (n100) [box, draw=black,  below of=n11] {\small Discrete $\mathcal{E}(k)$};
    \node (n101) [mybox, draw=black,  above left of=n100, xshift = - 1.15cm, yshift = -1cm] {\small Nesterov's \, Acceleration};
    \node (n110) [mybox, draw=black,  below left of=n100, xshift = - 1.15cm, yshift = +1cm] {\small Gradient Norm Minimization};

    \draw [arrow] (n00) --node [left] {dimensional analysis} (n10);
    \draw [arrow] (n10) -- (n11);
    \draw [arrow] (n11) -- node [right] {phase-space representation} (n100);
    \draw [arrow] (n100)  to [bend right=25] (n101);
    \draw [arrow] (n100)  to [bend left=25] (n110);
    \draw [arrow] [dashed] (n00) -- (n100);
    
\end{tikzpicture}
}
\end{center}
\caption{An illustration of our high-resolution ODE framework. The three solid straight lines represent Steps 1, 2 and 3, and the two curved lines denote Step 4. The dashed line is used to emphasize that it is difficult, if not impractical, to construct discrete Lyapunov functions directly from the algorithms.}
\label{fig:chart}
\end{figure}

\begin{figure}[htp!]
\centering
\includegraphics[scale=0.22]{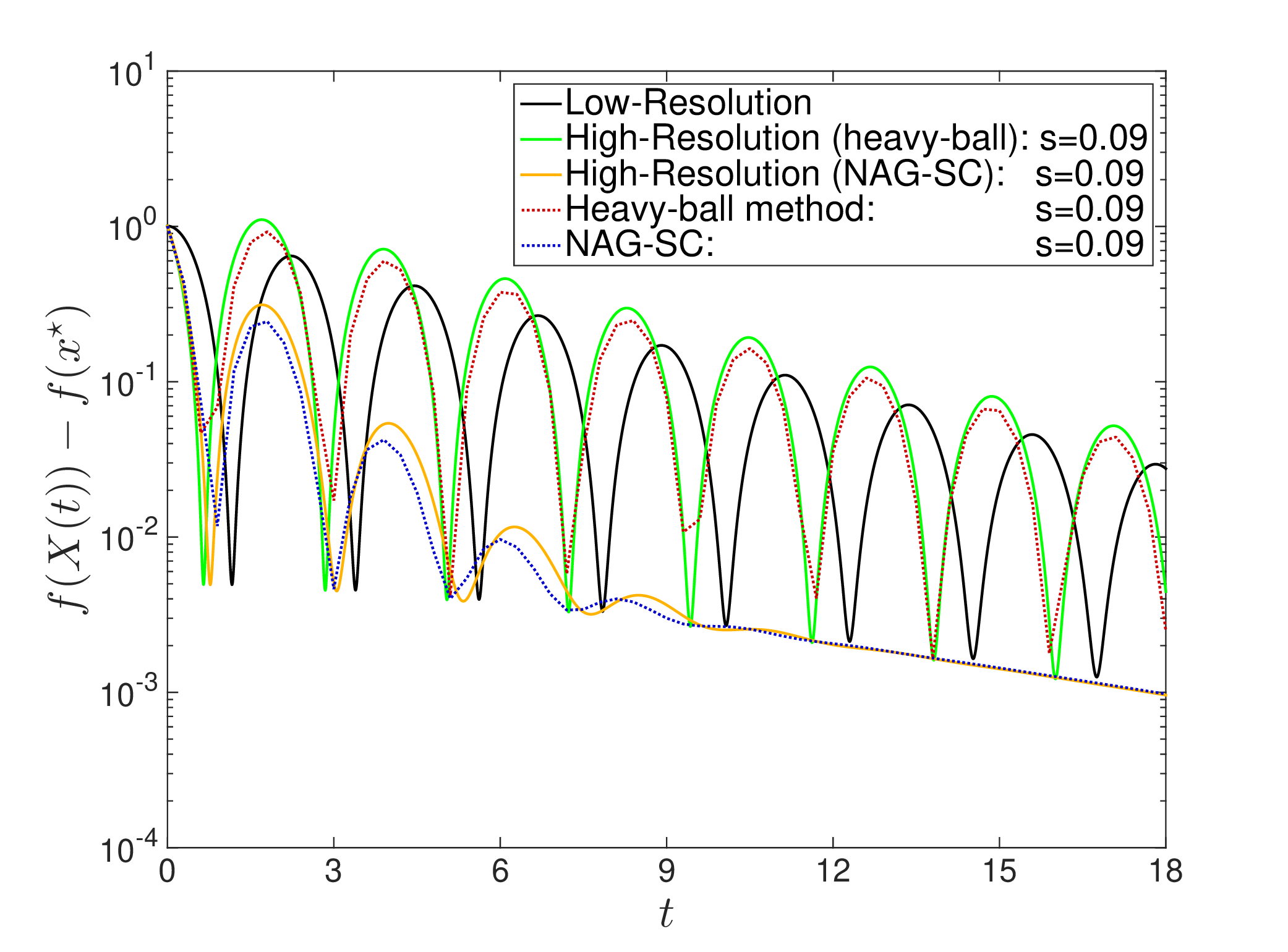}
\caption{The convergence behaviors of gradient-based algorithms and their approximating ODEs.} 
\label{fig: con-rate-sc}
\end{figure}

\paragraph{Step-1: Deriving high-resolution ODEs}  It can be observe that the second-order derivative can be approximated by the second-order central difference with a truncation error of $O(s)$, expressed as: 
\begin{equation}
\label{eqn: 2nd-order-approximation}
\frac{x_{k+1} + x_{k-1} - 2x_{k}}{s} = \ddot{X} + \text{O}(s).
\end{equation}
Thus, dimensional analysis from physics (as exemplified in~\cite{pedlosky2013geophysical}) was introduced in \cite{shi2022understanding} to derive high-resolution ODEs. For~\texttt{NAG-SC}~\eqref{eqn: nag-sc}, the high-resolution ODE is given by:
\begin{equation}
\label{eqn: nag-sc-high}
\ddot{X} + 2\sqrt{\mu} \dot{X} + \sqrt{s} \nabla^2 f(X) \dot{X} + (1 + \sqrt{\mu s}) \nabla f(X) = 0.
\end{equation}
Meanwhile, for the heavy-ball method~\eqref{eqn: hb}, the high-resolution ODE is
\begin{equation}
\label{eqn: hb-high}
\ddot{X} + 2\sqrt{\mu} \dot{X}  + (1 + \sqrt{\mu s}) \nabla f(X) = 0. 
\end{equation}
In contrast, the low-resolution ODE derived in~\cite{su2016differential} is identical for both the methods:
\begin{equation}
\label{eqn: nag-sc-low}
\ddot{X} + 2\sqrt{\mu} \dot{X} + \nabla f(X) = 0
\end{equation}
By comparing these three continuous-time models, we find that the high-resolution ODEs successfully differentiate between the two algorithms. This distinguishing capability is validated by the numerical demonstrations illustrated in~\Cref{fig: con-rate-sc}.

\paragraph{Step-2: Analyzing continuous dynamics via Lyapunov functions} Based on the high-resolution ODEs, we next construct appropriate Lyapunov functions to characterize their continuous-time dynamics. For the high-resolution ODE~\eqref{eqn: nag-sc-high}, we consider the Lyapunov function
\begin{align}
\mathcal{E}(t) = & \left( 1 + \sqrt{\mu s} \right) \left( f(X) - f(x^{\star}) \right) + \frac{1}{4} \| \dot{X} \|^2   \nonumber \\
                          & + \frac{1}{4}  \| \dot{X}  + 2\sqrt{\mu}(X - x^{\star}) + \sqrt{s} \nabla f(X) \|^2.            \label{eqn: lyapunov-ode-nag-sc}
\end{align}
Conversely, for the high-resolution ODE~\eqref{eqn: hb-high}, we define the Lyapunov function as
\begin{equation}
\label{eqn: lyapunov-ode-hb}
\mathcal{E}(t) = \left( 1 + \sqrt{\mu s} \right) \left( f(X) - f(x^{\star}) \right) + \frac{1}{4} \| \dot{X} \|^2 + \frac{1}{4} \| \dot{X} + 2\sqrt{\mu}(X - x^{\star}) \|^2.
\end{equation}
The underlying construction involves substituting the mixed term and its derivative into the high-resolution ODEs to isolate a term containing only the gradient. This is followed by integrating the kinetic energy and, ultimately, the potential energy. The details of this design principle can be found in~\cite{chen2025revisiting}. 

Both Lyapunov functions~\eqref{eqn: lyapunov-ode-nag-sc} and~\eqref{eqn: lyapunov-ode-hb} yield the same continuous-time convergence rate:
\[
f(X(t)) - f(x^{\star}) \leq O\left( e^{- \frac{\sqrt{\mu} t}{4}}\right). 
\]
However, a comparison reveals that the derivative inequality for the high-resolution ODE~\eqref{eqn: nag-sc-high} involves additional terms generated by the gradient correction, expressed as:
\begin{equation}
\label{eqn: deriv-nag-sc}
\frac{d \mathcal{E}}{dt} \leq - \frac{\sqrt{\mu}}{4} \mathcal{E} - \frac{\sqrt{s}}{2} \left[ \| \nabla f(X) \|^2 + \dot{X}^{\top} (\nabla^2 f) \dot{X} \right],
\end{equation}
where the extra structure demonstrates how the gradient correction term in the discretization can effectively accommodate a larger step size $s>0$.

\paragraph{Step-3: Constructing corresponding discrete Lyapunov functions} To transfer the continuous Lyapunov functions to the discrete counterparts, we introduce the velocity variable: $x_{k} - x_{k_1} = \sqrt{s} v_{k-1} $. By virtue of this substitution, the two gradient-based methods can be reformulated as the phase-space representation. Specifically,~\texttt{NAG-SC}~\eqref{eqn: nag-sc} can be rewritten as
\begin{equation}
    \label{eqn: hb-symplectic}
    \left\{ \begin{aligned}
            & x_{k} - x_{k-1} = \sqrt{s} v_{k-1}, \\
            & v_{k} - v_{k-1} = - \frac{2\sqrt{\mu s}}{1 - \sqrt{\mu s}} v_k - \frac{1 + \sqrt{\mu s}}{1 - \sqrt{\mu s}} \cdot \sqrt{s} \nabla f(x_k) \\
            & \mathrel{\phantom{v_{k} - v_{k-1} =}} - \sqrt{s} \left( \nabla f(x_k) - \nabla f(x_{k-1}) \right).
            \end{aligned} \right.
\end{equation}
Similarly, the heavy-ball method~\eqref{eqn: hb} can be expressed in terms of the velocity variable as the following discrete scheme:
\begin{equation}
    \label{eqn: nag-sc-symplectic}
    \left\{ \begin{aligned} 
             & x_{k} - x_{k-1} = \sqrt{s} v_{k-1}, \\
             & v_{k} - v_{k-1} = - \frac{2\sqrt{\mu s}}{1 - \sqrt{\mu s}} v_k - \frac{1 + \sqrt{\mu s}}{1 - \sqrt{\mu s}} \cdot \sqrt{s} \nabla f(x_k).
              \end{aligned} \right.
\end{equation}
In the absence of friction terms (i.e., setting the velocity coefficients to zero), the two discrete schemes above correspond exactly to the symplectic discretization of a Hamiltonian system \cite{feng1995collected}.  With the two phase-space representations, we can translate the two continuous Lyapunov functions,~\eqref{eqn: lyapunov-ode-nag-sc} and~\eqref{eqn: lyapunov-ode-hb}, into the discrete ones as: 
\begin{align}
\mathcal{E}(k) = & \frac{1 + \sqrt{\mu s}}{1 - \sqrt{\mu s}}  \left( f(x_k) - f(x^{\star}) \right) + \frac{1}{4} \left\| v_k \right\|^2 \nonumber \\ &+ \frac{1}{4} \left\| v_k + \frac{2\sqrt{\mu}}{1 - \sqrt{\mu s}}(y_{k+1} - x^{\star}) + \sqrt{s} \nabla f(y_k) \right\| - \underbrace{ \frac{s\| \nabla f(y_k) \|^2 }{2(1 - \sqrt{\mu s})}}_{\mathbf{additional\; term}}    \label{eqn: lyapunov-nag-sc}
\end{align}
and 
\begin{align}
\mathcal{E}(k) =&  \frac{1 + \sqrt{\mu s}}{1 - \sqrt{\mu s}}  \left( f(x_k) - f(x^{\star}) \right) + \frac{1}{4} \| v_k \|^2 \nonumber \\
                          & + \frac14  \left\| v_k + \frac{2\sqrt{\mu}}{1 - \sqrt{\mu s}} (x_{k+1} - x^{\star}) \right\|^2.   \label{eqn: lyapunov-hb}
\end{align}

\paragraph{Step-4: Establishing convergence rates through discrete Lyapunov functions}

Using the two discrete Lyapunov functions,~\eqref{eqn: lyapunov-nag-sc} and~\eqref{eqn: lyapunov-hb}, we can demonstrate that the permissible step size can be enlarged to $s = 1/(4L)$ for~\texttt{NAG-SC}, to $s = \mu/(16L^2)$ for the heavy-ball method~\eqref{eqn: hb}. Consequently, we derive the convergence rates for these two gradient-based methods to illustrate the acceleration mechanism, as detailed in the following theorems.

% is , thus, we derive the convergence rates for the two gradient-based methods to demonstrate the accelerated mechism, which are shown in the following two theorems. 

%\begin{table}[htp!]
%\centering
%\scalebox{0.8}{
% \begin{tabular}{|c|c|c|c|} 
%   \hline
%\diagbox{Objective}{Algorithms}& Gradient Descent &Polyak & Nesterov  \\ 
%  %\midrule 
%  \hline
%  $\mu$-strongly  & Global                                   & \textcolor{red}{Local}                                               &  \textcolor{red}{Global}  \\ 
%convex objective & $O\left( \left(1 - \frac{\mu}{L} \right)^k \right)$ & $O\left( \left(1 - \sqrt{\frac{\mu}{L} }\right)^k \right)$ &  $O\left( \left(1 - \sqrt{\frac{\mu}{L} }\right)^k \right)$ \\
%    \hline
% \end{tabular} 
% }
%\caption{The Convergence rate of gradient-based algorithms.} 
%\label{fig: convergence-rate}
%\end{table}

\begin{theorem}[Theorem 4 in \cite{shi2022understanding}]
\label{thm: heavy-ball}
Let $f \in \mathcal{S}_{\mu, L}^{1}$. If the step size is set to $s = \mu / (16L^2)$, the iterates $\{ x_k \}_{k=0}^{\infty}$ generated by the heavy-ball method satisfy
\[
f(x_k) - f(x_0) \leq \frac{5L \| x_0 - x^{\star} \|^2}{\left( 1 + \frac{\mu}{16L} \right)^k},
\]
for all $k \geq 0$. 
\end{theorem}

\begin{theorem}[Theorem 3 in \cite{shi2022understanding}]
\label{thm: nag-sc}
Let $f \in \mathcal{S}_{\mu, L}^{1}$. If the step size is set to $s = 1/(4L)$, the iterates $\{ x_k \}_{k=0}^{\infty}$ generates by~\texttt{NAG-SC} satisfy
\[
f(x_k) - f(x^{\star}) \leq \frac{5L \| x_0 - x^{\star} \|^2}{ \left( 1 + \frac{1}{12} \sqrt{ \frac{\mu}{L} } \right)^{k} },
\]
for all $k \geq 0$. 
\end{theorem}

In~\cite{chen2025revisiting}, the additional term cancels out, yielding a simplified form of the Lyapunov function~\eqref{eqn: lyapunov-nag-sc} for the implicit-velocity scheme~\eqref{eqn: nag-sc-y}: 
\begin{equation}
\label{eqn: lyapunov-nag-sc-1}
\mathcal{E}(k) = f(x_k) - f(x^{\star}) + \frac{\| v_k \|^2 }{4(1 + 2 \sqrt{\mu s} )} + \frac{\| v_k  + 2\sqrt{\mu} (x_k - x^{\star}) \|}{4},
\end{equation}
which leads directly to the following theorem.
\begin{theorem}[Theorem 4 in \cite{chen2025revisiting}]
\label{thm: nag-sc}
Let $f \in \mathcal{S}_{\mu, L}^{1}$. If the step size is set to $s = 1/L$, the iterates $\{ x_k \}_{k=0}^{\infty}$ generates by~\texttt{NAG-SC} satisfy
\[
f(x_k) - f(x^{\star}) \leq \frac{4L \| x_0 - x^{\star} \|^2}{ \left( 1 + \frac{1}{4} \sqrt{ \frac{\mu}{L} } \right)^{k} },
\]
for all $k \geq 0$. 
\end{theorem}

Subsequently, \cite{li2024linear} generalizes these results to composite optimization for proper, lower semi-continuous, and convex functions by deriving the following proximal inequality:
\[
\Phi(x - G_s(x)) \leq \Phi(y) + \left\langle G_s(x), x - y \right\rangle - \frac{s}{2} \| G_s(y) \|^2 - \frac{\mu}{2} \| y - x \|^2. 
\]

\section{NAG and FISTA}
\label{sec: nag-fista}

In this section, we consider the more general class of convex objective functions, denoted as $\mathcal{F}_L^{1}$, where the strong convexity assumption~\eqref{eqn: mu-strongly} is weaken to the standard convexity condition: for any $x, y \in \mathbb{R}^n$,
\begin{equation}
\label{eqn: convex}
f(y) \geq f(x) + \left\langle \nabla f(x), y - x \right\rangle,
\end{equation}
Consequently, the strong convexity terms on the left-hand sides of inequalities~\eqref{eqn: l-smooth-mu-strong-equiv-1} and~\eqref{eqn: l-smooth-mu-strong-equiv-2} reduce to zero. In this setting, classical gradient descent no longer converges linearly, but rather achieves a sub-linear convergence rate of
\[
f(x_k) - f(x^{\star}) \leq O\left(\frac1k\right).
\]
The celebrated~\texttt{NAG}, proposed by Nesterov~\cite{nesterov1983method}, accelerates this convergence rate to 
\[
f(x_k) - f(x^{\star}) \leq O\left(\frac{1}{k^2} \right).
\]
The iteration scheme for~\texttt{NAG} is given by:
\begin{subequations}
\label{eqn: nag}
\begin{align}
& x_{k+1} = y_{k} - s\nabla f(y_{k}),                                                \label{eqn: nag-gradient}          \\
& y_{k+1} = x_{k+1} + \frac{k}{k+r+1} (x_{k+1} - x_{k}),                 \label{eqn: nag-momentum} 
\end{align}    
\end{subequations}
where $s>0$ represents the step size, and $r \geq 2$ is a tuning parameter  (with $r = 2$ corresponding to Nesterov's original scheme). In the remainder of this section, we provide a brief review of several recent results analyzed through the lens of the high-resolution ODE framework~\cite{shi2022understanding}.

%=====================================================%
\subsection{Gradient norm minimization and composite optimization}
\label{subsec: composite}

In~\cite{shi2022understanding}, the high-resolution ODE for~\texttt{NAG}~\eqref{eqn: nag} is derived as
\begin{equation}
\label{eqn: high-res-ode}
\ddot{X} + \frac{r}{t} \dot{X} + \sqrt{s} \nabla f(X) \dot{X} + \left( 1 + \frac{3\sqrt{s}}{2t} \right) \nabla f(X) = 0. 
\end{equation}
Following the same four-step process outlined in~\Cref{sec: acceleration}, we construct the discrete Lyapunov function, using $v_{k-1} = (x_k - x_{k-1})/s$, as
\begin{multline}
\label{eqn: lyapunov-nag}
\mathcal{E}(k) = sk(k+r) \left( f(x_k) - f(x^{\star}) \right) \\+ \frac{1}{2} \| \sqrt{s} k v_{k-1} + r(y_{k} - x^{\star}) + sk \nabla f(y_{k-1}) \|^2
\end{multline}
which yields an accelerated convergence rate of $O(1/k^2)$ for the objective function error $f(x_k) - f(x^{\star})$. Furthermore, this framework allows us to establish a faster convergence rate for the gradient norm, achieving:
\begin{equation}
\label{eqn: nag-2k}
\min_{0 \leq i \leq k} \| \nabla f(x_i) \|^2 \leq O\left(\frac{1}{k^3}\right). 
\end{equation}
Unlike the objective function error, which requires knowing the optimal value $f(x_\star)$, the gradient norm $\|\nabla f(x_k)\|$ can be explicitly computed at each iteration. This makes it a highly practical quantity for monitoring convergence and establishing stopping criteria in real-world applications. 

\begin{figure}[htpb!]
\centering
\begin{minipage}[t]{0.45\linewidth}
\centering
\includegraphics[scale=0.38]{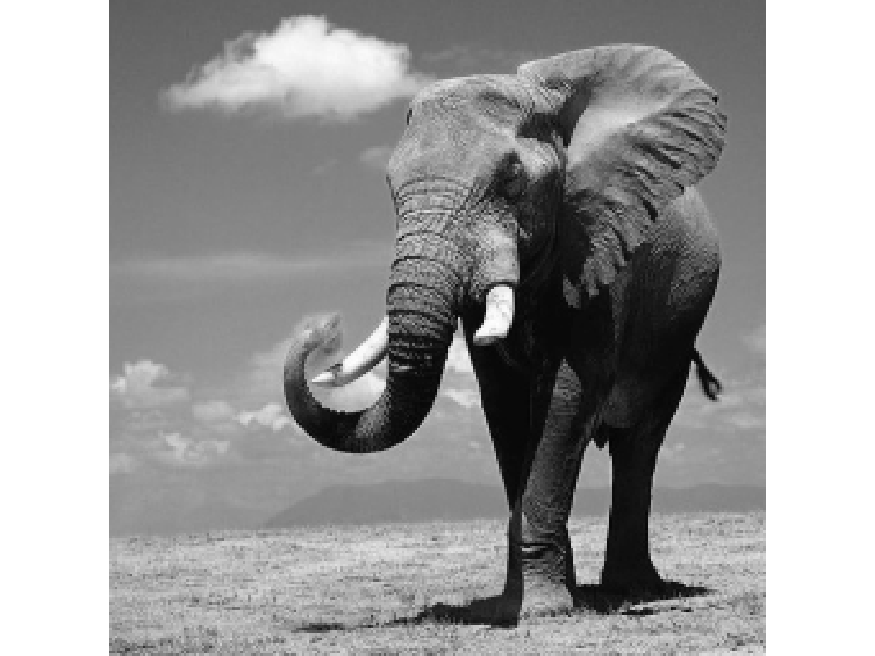}
\caption*{Original}
\end{minipage}
\begin{minipage}[t]{0.45\linewidth}
\centering
\includegraphics[scale=0.38]{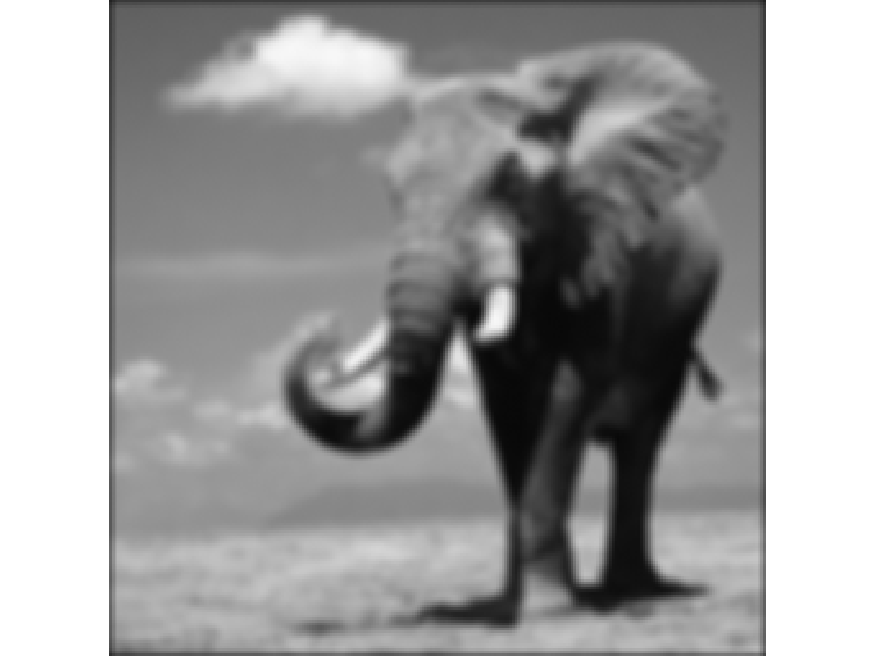}
\caption*{Blurred and Noisy}
\end{minipage}
%\begin{subfigure}[t]{0.325\linewidth}
%\centering
%\includegraphics[scale=0.3]{fig/recovery.eps}
%\caption{Deblurred and Denoisy}
%\end{subfigure}
%\caption{Deblurring and Denoising an image of an elephant by ISTA and FISTA.} 
\caption{The image deblurring problem: An image of an elephant}
\label{fig: deblurring}
\end{figure}

In \cite{chen2022gradient}, both the gradient-correction scheme and the implicit-velocity scheme for~\texttt{NAG}~\eqref{eqn: nag} was studied further, improving the convergence rate of the gradient norm to $o(1/k^3)$.  Equivalent forms of the discrete Lyapunov function~\eqref{eqn: lyapunov-nag} have been demonstrated for both schemes. For the implicit-velocity scheme with$v_k = (y_k - y_{k-1})/s$, , the Lyapunov function is expressed as:
\begin{equation}
\label{eqn: lyapunov-nag-1}
\mathcal{E}(k) = sk(k+r) \left( f(x_k) - f(x^{\star}) \right) + \frac{1}{2} \| \sqrt{s} (k - 1) v_{k} + r(x_{k} - x^{\star}) \|^2.
\end{equation}
Meanwhile, the original~\texttt{NAG} scheme~\eqref{eqn: nag}, it is given by:
\begin{equation}
\label{eqn: lyapunov-nag-2}
\mathcal{E}(k) = sk(k+r) \left( f(x_k) - f(x^{\star}) \right) + \frac{1}{2} \| \sqrt{s} k (y_k - x_k) + r(x_{k} - x^{\star}) \|^2.
\end{equation}
In \cite{li2026proximal}, the fundamental proximal inequality has been improved as
\begin{equation}
\label{eqn: proximal-convex}
\Phi(x - G_s(x)) \leq \Phi(y) + \left\langle G_s(x), x - y \right\rangle - \left( s - \frac{Ls^2}{2} \right)  \| G_s(y) \|^2. 
\end{equation} 
The accelerated gradient norm rate~\eqref{eqn: nag-2k} is generalized to composite optimization problems for~\texttt{FISTA}, an application of which is demonstrated in the image deblurring problem shown in~\Cref{fig: deblurring}. Letting $\mathcal{F}^{0}$ denote the class of proper, convex and lower semi-continuous functions, this result is rigorously stated in the following theorem.

\begin{theorem}[Theorem 5.1 in \cite{li2026proximal}]
\label{thm: fista}
Let $\Phi = f + g$ be a composite function, where $f \in \mathcal{F}_{L}^1$ and $g \in \mathcal{F}^0$. Then, for any step size $0 < s \leq 1/L$, the iterative sequence $\{x_k\}_{k=0}^{\infty}$ generated by \texttt{FISTA} satisfies
    \begin{equation}
    \label{eqn: fista-obj}
 \Phi(x_k) - \Phi(x^\star) \leq \frac{r^2\|x_0 - x^\star\|^2}{2sk(k+r)}.   
    \end{equation}
Moreover, if the step size satisfies $0< s< 1/L$, the iterative sequence $\{y_k\}_{k=0}^{\infty}$ satisfies
\begin{equation}
\label{eqn: fista-grad}
\left\{\begin{aligned}
& \min_{0\leq i \leq k} \| G_s(y_i) \|^2 \leq \frac{6r^2\|x_0 - x^\star\|^2}{s^2(1-sL) (k+1)\left(2k^2 + (6r+1)k + 6r^2 \right)},\\
& \lim_{k\rightarrow \infty} \left(k^3\min_{0\leq i \leq k} \| G_s(y_i) \|^2\right) = 0.
\end{aligned} \right.
\end{equation}
Furthermore, if $r>2$, the iterative sequence $\{x_k\}_{k=0}^{\infty}$ also satisfies
\begin{equation}
\label{eqn: fista-obj-fast}
\lim_{k\rightarrow \infty} \left[k^2 \left(\min_{0\leq i \leq k} \Phi(x_k) - \Phi(x^\star)\right)\right] = 0.
\end{equation}
\end{theorem}

%=====================================================%
\subsection{Underdamped case: $r \in [-1, 2)$}
\label{subsec: underdamped}

While the tuning parameter for the overdamped case ($r \geq 2$) has been thoroughly investigated via the high-resolution ODE framework in \cite{shi2022understanding, chen2022gradient, li2026proximal}, the underdamped case where $r \in [-1, 2)$ was analyzed in \cite{chen2023underdamped}.  Motivated by the intuition of scaling the $k-1$ factor in the mixed term to $(k-1)^\gamma$, we construct the discrete Lyapunov function as follows:
\begin{align}
\mathcal{E}(k) = & s^{\gamma}k^{\gamma} \left[k^{\gamma} + 2\gamma (k+1)^{\gamma}\right] \left( f(x_{k}) - f(x^{\star}) \right)                       \nonumber \\
                           & + \frac12 \left\| s^{\gamma/2}(k-1)^{\gamma}v_{k} + 2\gamma s^{(\gamma-1)/2}k^{\gamma-1}(x_k - x^{\star}) \right\|^2     \nonumber \\
                           & + \left\{  \frac{(k-2)k^{\gamma-1} \left[ (k-1)^{\gamma} + 2\gamma k^{\gamma-1} \right] }{k + 3\gamma - 2}  - (k-1)^{\gamma-1} (k-2)^{\gamma} \right\} \nonumber \\
                           & \mathrel{\phantom{+}} \cdot  2\gamma s^{\gamma-1} \left\| x_{k-1} - x^{\star} \right\|^2. \label{eqn: lyapunov-underdamped}
\end{align}
By analyzing this Lyapunov function, the convergence rate for the objective function error is derived as:
\begin{equation}
\label{eqn: fe-error-under}
f(x_k) - f(x^{\star}) \leq \text{O}\left( \frac{1}{k^{\frac{2(r+1)}{3}}} \right),
\end{equation}
and the convergence rate for the minimum gradient norm is given by:
\begin{equation}
\label{eqn: gn-error-under}
\min_{0 \leq i \leq k} \| \nabla f(x_i) \|^2 \leq \text{O}\left( \frac{1}{k^{\frac{2r+5}{3}}} \right). 
\end{equation}
Moreover, these two convergence rates,~\eqref{eqn: fe-error-under} and~\eqref{eqn: gn-error-under}, have been successfully generalized to the composite optimization for~\texttt{FISTA} in \cite{chen2023underdamped}.

%======================================================%
\subsection{Linear convergence on strongly convex functions}
\label{subsec: linear}

Take a closer look at the high-resolution ODE~\eqref{eqn: high-res-ode}, we observe that if the objective function is strongly convex, then the eigenvalues of the Hessian has a uniform lower bound: 
\[ 
| \lambda ( \nabla^2 f(x)) | \geq \mu.
\] 
Consequently, the friction does not diminish as time $t \rightarrow \infty$. If the step size satisfies $\sqrt{s} \mu \ll \sqrt{\mu}$, it is highly probable that linear convergence can be derived.

Based on the intuition from the high-resolution ODE~\eqref{eqn: high-res-ode}, a iteration-varying Lyapunov function is constructed in \cite{li2024linear2} as
\begin{multline}
\mathcal{E}(k) = \frac{s(k+r)(2k+r)}{1 - \mu s}  \left( f(x_k) - f(x^{\star}) \right) \\
+  \frac{s(k-1)^{2}}{2}  \left\|v_k \right\|^2 + \frac{1}{2} \left\| \sqrt{s}(k-1)v_k  + r(x_k - x^{\star}) \right\|^2, \label{eqn: lyapunov-c-sc}
\end{multline}
which establishes the linear convergence rate for~\texttt{NAG}~\eqref{eqn: nag}. Moreover, for the composite function $\Phi = f + g$ (where $f \in \mathcal{F}_{L}^1$ and $g \in \mathcal{F}^0$), the following proximal inequality is established in \cite{li2024linear2}: 
\begin{equation}
\| G_s(y)\|^2 \geq 2\mu \left( \Phi(y - s G_s(y)) - \Phi(x^{\star}) \right).  \label{eqn: proximal-c-sc}
\end{equation}
Thus, the linear convergence for~\texttt{FISTA} is established in \cite{li2024linear2}, resolving the open problem posed in \cite{chambolle2016introduction}.  The rigorous statement is provided below. 
\begin{theorem}[Theorem 4.4 in \cite{li2024linear2}]
\label{thm: fista-sc}
Let $\Phi = f + g$ be a composite function, where $f \in \mathcal{F}_{L}^1$ and $g \in \mathcal{F}^0$. For any step size $0 < s < 1/L$, there exists some positive integer $K:=K(L, \mu, s, r)$ such that the iterative sequence $\{ (x_k, y_k)\}^{\infty}_{k=0}$ generated by \texttt{FISTA} with any initial $x_0 = y_0 \in  \mathbb{R}$  satisfies the following inequalities as
\begin{equation}
\label{eqn: fista-sc}
\left\{ \begin{aligned}
         & \Phi(x_k) - \Phi(x^{\star})  \leq \frac{\mathcal{E}(K)}{s (k+r)(2k+r) \left[ 1 + (1- Ls) \cdot \frac{\mu s}{4} \right]^{k-K}}      \\
         & \| G_s(y_k)\|                     \leq \frac{4\mathcal{E}(K)}{s^{2}(1 - Ls) (k+r)(2k+r) \left[ 1 + (1- Ls) \cdot \frac{\mu s}{4} \right]^{k-K}}   
         \end{aligned} \right. 
\end{equation}
for any $k \geq K$.
\end{theorem}

%======================================================%
\subsection{Monotonicity}
\label{subsec: monotonicity}

\begin{figure}[htbp!]
\centering
\scalebox{0.8}{
\begin{tikzpicture}[
  node distance=1.2cm and 1.5cm,
  box/.style={draw, minimum width=2.4cm, minimum height=1.0cm, align=center},
  smallbox/.style={draw, minimum width=4.2cm, minimum height=1.8cm, align=center},
  every node/.style={font=\sffamily}
]

\node[smallbox, transform shape, scale=0.8] (momentum) {
  Momentum Update~\eqref{eqn: nag-momentum} \\[2pt] 
  $ y_{k+1} = x_{k+1} + \dfrac{k}{k+r+1} (x_{k+1} - x_{k}) $
};

\node[smallbox, below=3.6cm  of momentum, transform shape, scale=0.8] (gradient) {
  Gradient Update~\eqref{eqn: nag-gradient} \\[2pt] 
  $x_{k+1} = y_{k} - s \nabla f(y_k)$
};

\node[smallbox, below=0.8cm of gradient, transform shape, scale=0.8] (mnag) {\texttt{M-NAG}~\eqref{eqn: essential-iteration}\\[2pt] $y_{k+1}  = x_{k+1} +  \frac{k}{k+r+1}(x_{k+1} - x_{k}) + \frac{k+r}{k+r+1}(y_k - s \nabla f(y_k) - x_{k+1})$};
\node[box, right =2.8cm of momentum, yshift=-2.6cm, transform shape, scale=0.8] (nag) {\texttt{NAG}~\eqref{eqn: nag}};

\node[draw, circle, minimum size=1.2cm, transform shape, scale=0.8] (fraction) at ($(momentum)!0.55!(gradient) $) {\( \pmb{\dfrac{k + r}{k + r + 1}}\)};
\node at ($(fraction.north)!0.4!(momentum.south)$) {\Large$+$};  
\node at ($(fraction.south)!0.5!(gradient.north)$) {\Large$\times$};  % Near circle-gradient
\node at ($(fraction.north)!1.08!(gradient.south)$) {\rotatebox{90}{\Large$=$}};

% Momentum to NAG: right then down
\coordinate (momentum-east) at ([xshift=3.75cm]momentum.east);
\draw[thick, -{Stealth}] (nag)-- (momentum-east) -- (momentum) ;

% Gradient to NAG: right then up
\coordinate (gradient-east) at ([xshift=4.65cm]gradient.east);
\draw[thick, -{Stealth}] (nag) -- (gradient-east) -- (gradient);

\end{tikzpicture}
}
\caption{Illustration of how~\texttt{M-NAG} consolidates the momentum and gradient steps of classical~\texttt{NAG} into a unified update rule. Each solid arrow $A \rightarrow B$ indicates that $B$ is a component or derived step of $A$.}
\label{fig: digram-nag-mnag}
\end{figure}
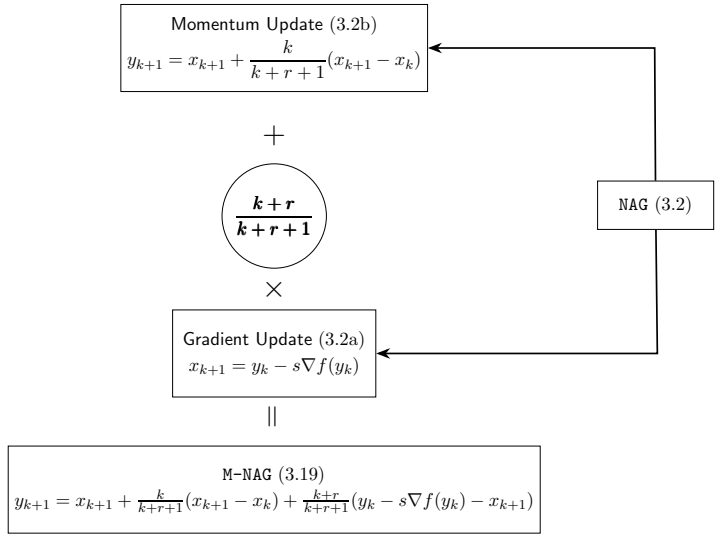

A monotonic variant of~\texttt{NAG}, referred to as~\texttt{M-NAG}, was proposed in~\cite{beck2009fast} and is defined by the following iteration scheme:
\begin{subequations}
\label{eqn: m-nag}
\begin{align}
& z_{k}       = y_{k} - s\nabla f(y_{k}),                                                                                                                              \label{eqn: m-nag-gradient} \\
& x_{k+1}   = \left\{ \begin{aligned} 
                               & z_{k}, && \text{if}\; f(z_{k}) \leq f(x_{k}), \\
                               & x_{k}, && \text{if}\; f(z_{k}) > f(x_{k}),
                               \end{aligned} \right.                                                                                                                          \label{eqn: m-nag-monetone} \\
& y_{k+1}  = x_{k+1} + \frac{k}{k+r+1} (x_{k+1} - x_{k}) + \frac{k+r}{k+r+1} (z_{k} - x_{k+1}).                                        \label{eqn: m-nag-momentum} 
\end{align}    
\end{subequations}
In the final part of this section, we review the theoretical insights established in \cite{fu2024lyapunov} based on the high-resolution ODE framework. By substituting~\eqref{eqn: m-nag-gradient} into~\eqref{eqn: m-nag-momentum}, we can express the relationship between the iterative sequences~\texttt{M-NAG}, $\{x_k\}_{k=0}^{\infty}$ and $\{y_{k}\}_{k=0}^{\infty}$ as follows: 
\begin{equation}
\label{eqn: essential-iteration}
y_{k+1}  = x_{k+1} +  \frac{k}{k+r+1}(x_{k+1} - x_{k}) + \frac{k+r}{k+r+1}(y_k - s \nabla f(y_k) - x_{k+1}),
\end{equation}
which shows that~\texttt{M-NAG} executes a single update that couples information from both sequences $\{x_k\}_{k=0}^{\infty}$ and  $\{y_k\}_{k=0}^{\infty}$. Crucially, this formulation shows how updating $\{y_k\}_{k=0}^{\infty}$ guarantees that the objective values $\{f(x_k)\}_{k=0}^{\infty}$ form a monotonically non-increasing sequence. Furthermore, this relation~\eqref{eqn: essential-iteration} reveals that~\texttt{M-NAG} is essentially a linear combination of the momentum step~\eqref{eqn: nag-momentum} and the gradient step~\eqref{eqn: nag-gradient}, as depicted in~\Cref{fig: digram-nag-mnag}. In this manner,~\texttt{M-NAG} retains only part information from the full update ~\texttt{NAG}~\eqref{eqn: nag}, effectively isolating and recombining its key components.

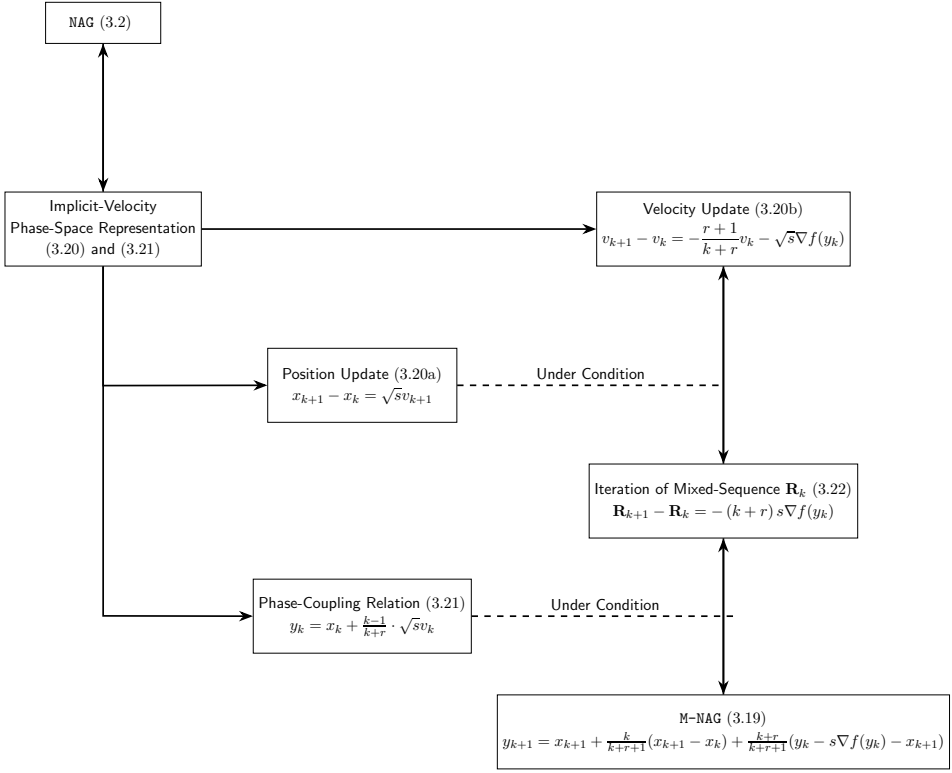
\begin{figure}[htbp!]
\centering
\scalebox{0.8}{
\begin{tikzpicture}[
  node distance=1.2cm and 1.5cm,
  box/.style={draw, minimum width=2.8cm, minimum height=1.0cm, align=center},
  smallbox/.style={draw, minimum width=4.6cm, minimum height=1.8cm, align=center},
  every node/.style={font=\sffamily},
  scale=0.68,
  transform shape
]

\node[box] (nag) {
  \texttt{NAG}~\eqref{eqn: nag}
};

% Phase-space representation
\node[smallbox, below=3.6cm of nag] (nag-phase) {
 Implicit-Velocity  \\[2pt] Phase-Space Representation \\[2pt] \eqref{eqn: nag-phase} and~\eqref{eqn: nag-phase-x-y-v} 
};

% Middle row
\node[smallbox, right=1.6cm of nag-phase, yshift=-3.8cm] (position) {
 Position Update~\eqref{eqn: nag-phase-x} \\[2pt]
  $ x_{k+1} - x_k = \sqrt{s} v_{k+1}  $
};

\node[smallbox, below=3.8cm of position] (coupling) {
  Phase-Coupling Relation~\eqref{eqn: nag-phase-x-y-v} \\[2pt]
  $ y_{k} = x_{k} + \frac{k-1}{k+r} \cdot \sqrt{s} v_{k} $
};

\node[smallbox, right=1.6cm of nag-phase, xshift=8cm] (velocity) {
   Velocity Update~\eqref{eqn: nag-phase-v} \\[2pt]
  $ v_{k+1} - v_{k} = - \dfrac{r+1}{k+r} v_{k} - \sqrt{s} \nabla f(y_k)$
};

% R_{k+1}
\node[smallbox, below=4.8cm of velocity] (mix) {
Iteration of Mixed-Sequence $\mathbf{R}_k$~\eqref{eqn: nag-iv-iter}\\[2pt] 
 $\mathbf{R}_{k+1} - \mathbf{R}_{k} = - \left( k+r\right)s \nabla f(y_{k})  $
};

% M-NAG
\node[smallbox, below=3.8cm of mix] (mnag) {
 \texttt{M-NAG}~\eqref{eqn: essential-iteration} \\[2pt] 
 $y_{k+1}  = x_{k+1} +  \frac{k}{k+r+1}(x_{k+1} - x_{k}) + \frac{k+r}{k+r+1}(y_k - s \nabla f(y_k) - x_{k+1})$
};

Arrows
\draw[thick, - {Stealth}] (nag) --  (nag-phase);
\draw[thick, -{Stealth}] (nag-phase) -- (nag);
\draw[thick, - {Stealth}] (velocity) --  (mix);
\draw[thick, -{Stealth}] (mix) -- (velocity);
\draw[thick, - {Stealth}] (mnag) --  (mix);
\draw[thick, -{Stealth}] (mix) -- (mnag);

% nag-phase to velocity (unchanged, horizontal)
\draw[thick, -{Stealth}] (nag-phase) -- (velocity);

% nag-phase to position:  vertical then horizontal 
\coordinate (nag-to-pos-x) at ([yshift=-2.9cm]nag-phase.south);
\draw[thick, -{Stealth}] (nag-phase.south) -- (nag-to-pos-x) -- (position.west);

% nag-phase to coulpling:  vertical then horizontal 
\coordinate (nag-to-coupling-x) at ([yshift=-8.5cm]nag-phase.south);
\draw[thick, -{Stealth}] (nag-phase.south) -- (nag-to-coupling-x) -- (coupling.west);

% position equivalent between velocity and mix
\coordinate (position-to-vm-x) at ([xshift=6.5cm]position.east);
\draw[thick, -, dashed] (position.east) -- node[midway, above]{Under Condition}(position-to-vm-x) ;

% coupling equivalent between mix and mnag
\coordinate (coupling-to-mm-x) at ([xshift=6.5cm]coupling.east);
\draw[thick, -, dashed] (coupling.east) -- node[midway, above]{Under Condition}(coupling-to-mm-x);

\end{tikzpicture}
}
\caption{Illustration of how~\texttt{M-NAG} recovers the~\texttt{NAG} update though the implicit-velocity phase-representation. Each solid arrow $A \rightarrow B$ indicates that $B$ is a component or derived step of $A$, while each dashed arrow represents a conditional relation required to establish the connection. }
\label{fig: digram-nag-phase-mnag}
\end{figure}

Under this formulation, the classical~\texttt{NAG} iteration can be rewritten as:
\begin{subequations}
\label{eqn: nag-phase}
\begin{align}
         & x_{k+1} - x_{k} = \sqrt{s} v_{k+1},                                                                               \label{eqn: nag-phase-x}          \\
         & v_{k+1} - v_{k} = -\frac{r+1}{k+r}  v_{k} - \sqrt{s} \nabla f\left( y_{k} \right),                 \label{eqn: nag-phase-v}
\end{align}    
\end{subequations}
where the auxiliary sequence $\{y_{k}\}_{k=0}^{\infty}$ satisfies the phase-coupling relation:
\begin{equation}
\label{eqn: nag-phase-x-y-v}
y_k = x_{k} +  \frac{k-1}{k+r} \cdot \sqrt{s}v_k.
\end{equation}
Following~\cite{chen2022gradient}, we introduce the mixed sequence $\mathbf{R}_k = (k-1)\sqrt{s}v_{k} + rx_k$. Combining~\eqref{eqn: nag-phase-x} with~\eqref{eqn: nag-phase-v} yields the update rule:
\begin{equation}
\mathbf{R}_{k+1} - \mathbf{R}_{k} = - \left( k+r\right)s \nabla f(y_{k})   \label{eqn: nag-iv-iter}. 
\end{equation}
The implicit-velocity phase-space formulation clarifies the precise connection between~\texttt{M-NAG} and~\texttt{NAG}.  As illustrated in~\Cref{fig: digram-nag-phase-mnag}, a bidirectional equivalence exists among the velocity update~\eqref{eqn: nag-phase-v}, the mixed-sequence iteration~\eqref{eqn: nag-iv-iter}, and the~\texttt{M-NAG} scheme~\eqref{eqn: essential-iteration}. Specifically:
\begin{itemize}
\item Given the position update~\eqref{eqn: nag-phase-x}, the velocity update~\eqref{eqn: nag-phase-v} is equivalent to the mixed-sequence iteration~\eqref{eqn: nag-iv-iter}. 
\item Given the phase-coupling relation~\eqref{eqn: nag-phase-x-y-v}, the mixed-sequence iteration~\eqref{eqn: nag-iv-iter} is equivalent to the~\texttt{M-NAG} update rule~\eqref{eqn: essential-iteration}. 
\end{itemize}
Since the velocity sequence $\{v_k\}_{k=0}^{\infty}$ is defined directly via $\{x_k\}_{k=0}^{\infty}$ and $\{y_k\}_{k=0}^{\infty}$, each of the two conditions,~\eqref{eqn: nag-phase-x} and~\eqref{eqn: nag-phase-x-y-v}, can be imposed independently without external assumptions. Crucially, however, enforcing both conditions simultaneously does not inherently follow from \texttt{M-NAG}.  Therefore, while \texttt{M-NAG} update~\eqref{eqn: essential-iteration} can recover the full NAG iteration, it requires an additional external assumption: either the position update~\eqref{eqn: nag-phase-x} or the phase-coupling relation~\eqref{eqn: nag-phase-x-y-v}. This highlights that \texttt{M-NAG} only partially encodes the full structure of \texttt{NAG} derived from the implicit-velocity formulation.

As shown above, the mixed-sequece iteration is logically equivalent to the~\texttt{M-NAG} update rule. By substituting the phase-coupling relation~\eqref{eqn: nag-phase-x-y-v} to the definition of $\mathbf{S}_k$, we eliminate the velocity term $v_{k}$, yielding:
\begin{align}
\mathbf{S}_k  & =  (k-1)\sqrt{s}v_{k} + r x_k - \left( k+r\right)s \nabla f(y_{k}) \nonumber \\ & = (k + r)y_{k}  - kx_k  - (k + r)s\nabla f(y_k). \label{eqn: mix-new}
\end{align}
Using this expression, the new iterative difference can be reformulated as:
\begin{multline}
(k + r + 1) y_{k+1}  - (k + 1) x_{k+1}  - \left( k + r + 1 \right)\nabla f(y_{k+1}) \\
= (k + r)y_{k}  - kx_k  - (k + r)s\nabla f(y_k) - \left( k + r + 1 \right)\nabla f(y_{k+1}), \label{eqn: m-nag-iteration-new} 
\end{multline}
which is exactly equivalent to the~\texttt{M-NAG} scheme~\eqref{eqn: essential-iteration}. Consequently,  the reformulated iteration~\eqref{eqn: m-nag-iteration-new} is adopted as the foundation for constructing the Lyapunov function in \cite{fu2024lyapunov}. Together with the earlier formulation~\eqref{eqn: nag-iv-iter}, the result demonstrate that the construction of a Lyapunov function fundamentally depends only on the \texttt{M-NAG} formulation, without requiring the full update structure of classical~\texttt{NAG}. Furthermore, according to the proximal inequality~\eqref{eqn: proximal-c-sc}, the smooth results can be generalized to the composite optimization via the following Lyapunov function
\begin{align}
\mathcal{E}(k) = &\; s(k+1)(k+r+1)\left( \Phi(x_{k+1}) - \Phi(x^{\star}) \right) \nonumber \\ 
                           &\; + \frac12\left\| (k-1)\sqrt{s}v_{k} + r (x_{k} - x^{\star}) - s (k + r) G_s(y_k) \right\|^2, \label{eqn: lyapunov-fista}
\end{align}
which establishes the linear convergence rate. The rigorous statement is presented as follows.

\begin{theorem}[Theorem 4.2 in \cite{fu2024lyapunov}]
\label{thm: fista-mono}
Let $\Phi = f + g$, where $f \in \mathcal{S}_{\mu,L}^{1}(\mathbb{R}^d)$ and $g \in \mathcal{F}^0(\mathbb{R}^d)$. Given that any step size satisfies $0 < s < 1/L$, there exists a positive integer $K: = \max\left\{0,  \frac{3r^2 - 4r - 12}{8}\right\}$ such that the iterative sequence $\{x_{k}\}_{k=0}^{\infty}$ generated by~\texttt{M-FISTA},  with any initial $x_0 = y_0 \in \mathbb{R}^d$, satisfies the following inequality: 
\begin{equation}
\label{eqn: fista-rate}
\Phi(x_k) - \Phi(x^{\star}) \leq \frac{(r +1)\left( \Phi(x_1) - \Phi(x^{\star})\right) + r^2 L\| x_1 - x^{\star} \|^2}{k(k+r) \left[ 1 + (1 - Ls) \cdot \frac{\mu s}{4} \right]^k}, 
\end{equation}
for any $k \geq \max\left\{ 1, K \right\}$. 
\end{theorem}

Moreover, the results regarding the critical step size $s = 1/L$ within the high-resolution ODE framework are established in~\cite{fu2025family}, which also suggests a potentially faster convergence rate in this direction.

\section{ADMM and PDHG}
\label{sec: admm-pdhg}

In this section, we review recent advances in high-resolution ODE frameworks designed for analyzing advanced gradient-based algorithms, with a focus on \texttt{ADMM} and \texttt{PDHG}. In modern applications such as total-variation denoising and $\ell_1$ trend filtering, the generalized least absolute shrinkage and selection operator (\texttt{Lasso}) is frequently formulated as:
\[
\min_{x \in \mathbb{R}^d} \Phi(x):=\|Ax - b\|^2 + \lambda\|Fx\|_1,
\]
where the standard gradient-based or proximal gradient algorithms cannot directly solve this problem. Because the linear operator $F$ is nested inside the non-smooth regularization term, an efficient computation of its proximal operator is prevented. More generally, this problem can be framed in the following composite optimization form:
\begin{equation}
\label{eqn: composite-gen-lasso}
\min_{x \in \mathbb{R}^d} \Phi(x):= f(x) + g(Fx),
\end{equation}
where $f$ is smooth and convex, and $g$ is a proper, convex, and lower semi-continuous function. 

%=================================%
\subsection{ADMM and its mechanism}
\label{subsec: admm-mechanism}

To decouple the linear operator from the non-smooth term, we can express~\eqref{eqn: composite-gen-lasso} in the general two-block constrained form:
\begin{equation}
\label{eqn: admm-problem}
  \min \; f(x) + g(y) \qquad \mathrm{s.t.}     \quad  Fx + Gy = h
\end{equation}
where $x \in \mathbb{R}^{d_1}$ and $y \in \mathbb{R}^{d_2}$ are primal variables, $F \in \mathbb{R}^{m \times d_1}$ and $G \in \mathbb{R}^{m \times d_2}$ are two matrices, and $h \in \mathbb{R}^{m}$ is a $m$-dimensional vector. The~\texttt{ADMM} iterations are derived from the augmented Lagrangian:
\begin{equation}
\label{eqn: augment-lagrangian}
L_{s}(x,y; \lambda) = f(x) + g(y) + \left\langle \lambda, Fx + Gy -h \right\rangle + \frac{1}{2s} \left\| Fx + Gy -h  \right\|^2 ,
\end{equation}
where $\lambda \in \mathbb{R}^m$ is the dual multiplier and $s>0$ is a penalty parameter. The algorithm proceeds by sequentially updating $x$, $y$, and the dual multiplier $\lambda$ according to the following steps:    

\begin{subequations}
\label{eqn: admm}
\begin{align}
         & x_{k+1} = \argmin_{x \in \mathbb{R}^{d_1}} \left\{ f(x) + \frac{1}{2s} \|Fx + Gy_k - h + s \lambda_k\|^2  \right\},                   \label{eqn: admm-x}            \\   
         & y_{k+1} = \argmin_{y \in \mathbb{R}^{d_2}} \left\{ g(y) + \frac{1}{2s} \|Fx_{k+1} +  Gy - h + s \lambda_k\|^2 \right\},            \label{eqn: admm-y}            \\
         & \lambda_{k+1} = \lambda_k + \frac1s (Fx_{k+1} + Gy_{k+1} - h).                                                                                          \label{eqn: admm-lambda}
\end{align}    
\end{subequations}

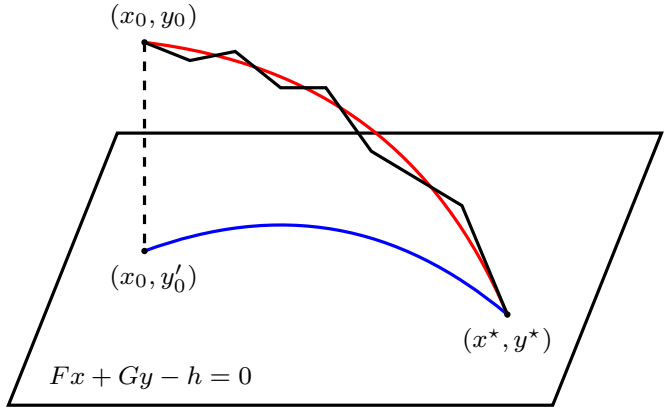
\begin{figure}[htb!]
  \centering
  \scalebox{1.2}{\begin{tikzpicture}[line width=1.0pt].
    % Plane
    \draw[xslant=0.4](-1,0) rectangle (5,3);
    \node at(-0.7,0.3)[right]{\scriptsize\( Fx + Gy - h = 0 \)};
  
    % Low-Resolution ODE
    \path[blue] (0.5, 1.7) edge [bend left] (4.5,1);
    \fill(0.5,1.7)circle(1pt)node[below, xshift=0.1cm]{\scriptsize\( (x_0, y_0') \)};
%    \draw[-latex,thick](1.7,-0.2)node[below]{\scriptsize\( L-R\ ODE \)}--(2.7,1.2);

    % High-Resolution ODE
    \path[red] (0.5, 4) edge [bend left] (4.5,1);
%    \draw[-latex,thick](3.5,4.7)node[above]{\scriptsize\( H-R\ ODE \)}--(2.3,3.7);

    % ADMM
    \draw[black] (0.5,4) -- (1.0, 3.8) -- (1.5,3.9)-- (2.0,3.5) -- (2.5,3.5) -- (3.0, 2.8) -- (3.5,2.5) -- (4,2.2) -- (4.5,1);
    \fill(0.5,4)circle(1pt)node[above, xshift=0.1cm]{\scriptsize\( (x_0, y_0) \)};
%    \draw[-latex,thick](5.5,3.7)node[above,xshift=0.3cm]{\scriptsize\( ADMM \)}--(3.4,3.1);

    \fill(4.5,1)circle(1pt)node[below]{\scriptsize\( (x^\star, y^\star) \)};

    \draw[black,dashed](0.5,4) -- (0.5,1.7); 
  \end{tikzpicture}}
\caption{Schematic comparison of the discrete~\texttt{ADMM} trajectory (black), the continuous limit ODE (blue), and the high-resolution ODE system (red).  Unlike the low-resolution ODE, which requires the initialization $y_0'$ to satisfy $Fx_0 +Gy_0' = h$, the high-resolution model captures the trajectory for any arbitrary initial $(x_0, y_0)$, mirroring the discrete algorithm's flexibility.}
\label{fig: admm-low-high}
\end{figure}

As established in \cite{li2024understanding1}, the high-resolution ODE system modeling~\texttt{ADMM} is given by
\begin{subequations}
\label{eqn: admm-high-ode}
\begin{align}
        & F^{\top}G\dot{Y}  = F^{\top}\Lambda + \nabla f(X),             \label{eqn: admm-ode-x}               \\
        & 0 = G^{\top}\Lambda  + \nabla g(Y),                                   \label{eqn: admm-ode-y}               \\
        & s^2 \dot{\Lambda} = FX + GY - h.                                      \label{eqn: admm-ode-lambda}  
\end{align}    
\end{subequations}
Compared to the previous low-resolution ODE model, the primary advancement lies in~\eqref{eqn: admm-ode-lambda}, which explicitly captures the $\lambda$-correction effect, a distinctive feature of \texttt{ADMM}. The small but essential perturbation causes the variable pair $(X,Y)$ to deviate from the constraint hyperplane $Fx+ Gy= h$. This underlying dynamical phenomenon is visualized in~\Cref{fig: admm-low-high}.

The continuous dynamical model was established, which enables the transition to a discrete Lyapunov function:
\begin{equation}
\label{eqn: lyapunov}
\mathcal{E}(k) = \frac{1}{2s} \| G(y_k - y^{\star}) \|^2 + \frac{s}{2} \| \lambda_k - \lambda^{\star} \|^2,
\end{equation}
which leads to the convergence rate for both the primal-dual gap and the numerical error induced by implicit discretization. Furthermore, by defining the numerical error itself as a Lyapunov function: 
\begin{equation}
\label{eqn: lyapunov-ne}
\mathcal{NE}(k) = \frac{1}{2s} \| G(y_{k+1} - y_{k}) \|^2 + \frac{s}{2} \| \lambda_{k+1} - \lambda_{k} \|^2 
\end{equation}
which determines the numerical error is monotonically non-increasing. The rigorous statement  is established in \cite{li2024understanding1}, and is presented below.

\begin{theorem}[Theorem 5.1 and Theorem 5.3 in \cite{li2024understanding1}]
\label{thm: admm-convegence-rate}
Let $f \in \mathcal{F}^0(\mathbb{R}^d)$ and  $g \in \mathcal{F}^0(\mathbb{R}^d)$. For any step size $s>0$, the ergodic sequence $\{(\overline{x}_N, \overline{y}_N; \overline{\lambda}_N)\}_{N=0}^{\infty}$ generated by the~\texttt{ADMM} iterations~\eqref{eqn: admm-x} --- \eqref{eqn: admm-lambda} converges to the saddle point $(x^{\star},y^{\star};\lambda^{\star})$ with the following rate: 
\begin{equation}
\label{eqn: average-rate-admm}
\mathcal{L}(\overline{x}_{N}, \overline{y}_{N}; \lambda^{\star}) - \mathcal{L}(x^{\star}, y^{\star}; \overline{\lambda}_{N}) \leq \frac{\left\|G(y_0 - y^{\star})\right\|^2 + s^2 \left\| \lambda_0 - \lambda^{\star} \right\|^2}{2s(N+1)}.
\end{equation}
Moreover, for any initial $(x_0, y_0; \lambda_0) \in \mathbb{R}^{d_1} \times \mathbb{R}^{d_2} \times \mathbb{R}^{m}$, the iterative sequence $\{ (x_k, y_k; \lambda_k) \}_{k=0}^{\infty}$ satisfies  the following rate for the last iterate: 
\begin{align}
\left\|G(y_{N+1} - y_{N})\right\|^2 & + s^2 \left\| \lambda_{N+1} - \lambda_{N} \right\|^2 \nonumber \\ & \leq \frac{\left\|G(y_0 - y^{\star})\right\|^2 + s^2 \left\| \lambda_0 - \lambda^{\star} \right\|^2}{N+1}.  \label{eqn:admm-convegence-rate-mono}
\end{align}
\end{theorem}

%=================================%
\subsection{PDHG and its variants}
\label{subsec: pdhg-variant}
\begin{figure}[htb!]
\centering
\includegraphics[scale=0.20]{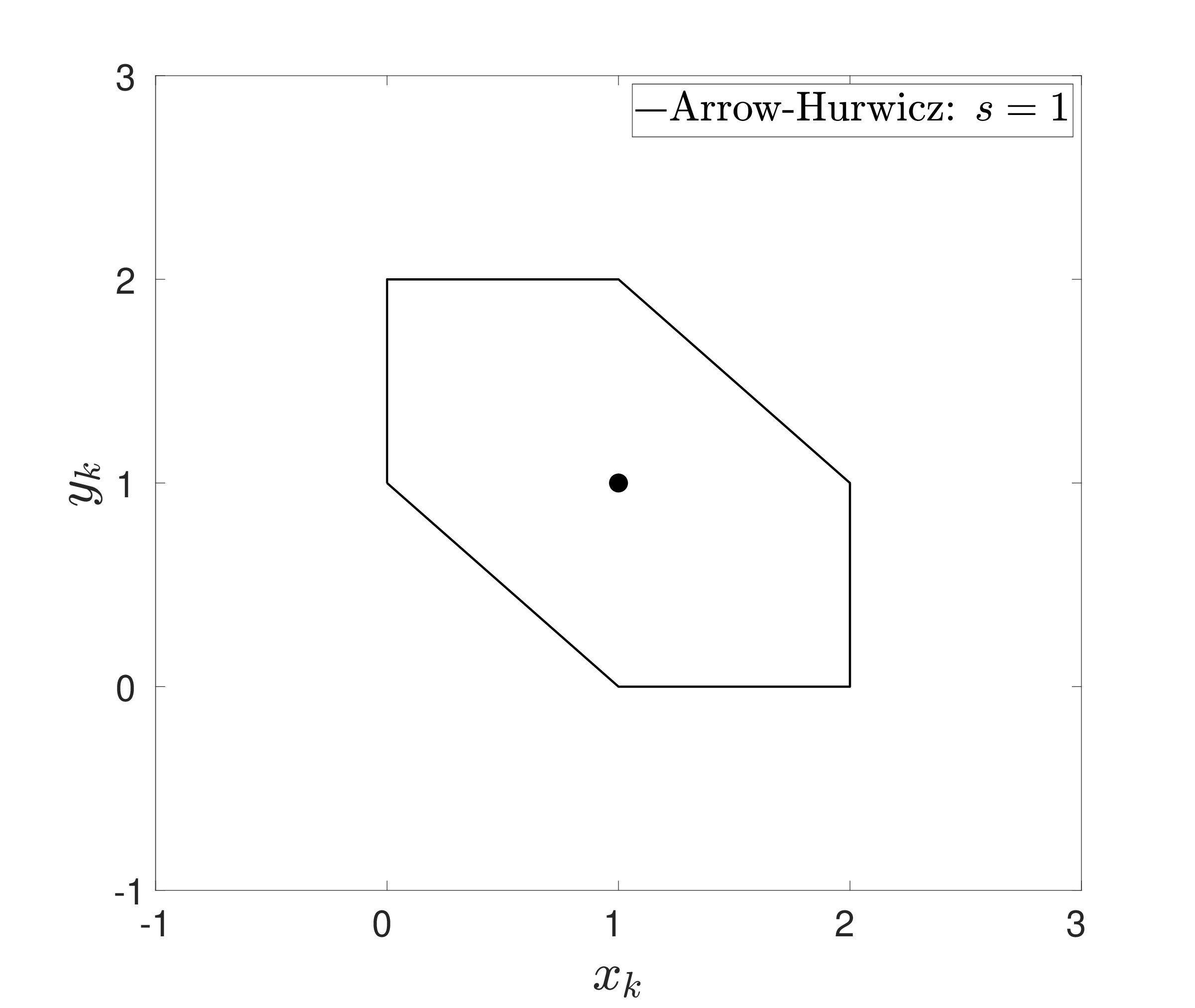}
\caption{A counterexample demonstrating the non-convergence of the Arrow-Hurwicz algorithm. For a bilinear saddle-point objective $\Phi(x, y) = x - xy + y$ and the initial point $(0, 1)$, the algorithm's iterates fail to converge to the unique saddle point $(1,1)$. }
\label{fig: s1}
\end{figure}

By utilizing the Legendre-Fenchel transform, the composite optimization problem~\eqref{eqn: composite-gen-lasso} can also be transformed into a minimax problem as:
\begin{equation}
\label{eqn: game-rep}
\min_{y \in \mathbb{R}^{d_2}} \max_{x \in \mathbb{R}^{d_1}} \Phi(x,y) = \min_{x \in \mathbb{R}^{d_1}} \max_{y \in \mathbb{R}^{d_2}} \Phi(x,y) := f(x) + \big\langle Fx, y \big\rangle - g^{\star}(y),
\end{equation}
which was original proposed in \cite{arrow1958studies}, alongside a saddle-point solving method known as the Arrow–Hurwicz algorithm:
\begin{subequations}
\label{eqn: arrow-hurwicz}
\begin{align}
  & x_{k+1} = \argmin_{x \in \mathbb{R}^{d_1}} \left\{ f(x) + \big\langle Fx, y_k \big\rangle + \frac{1}{2s} \|x - x_k \|^2 \right\},                         \label{eqn: ah-descent} \\
  & y_{k+1} = \argmax_{y \in \mathbb{R}^{d_2}} \left\{ - g^{\star}(y) + \big\langle Fx_{k+1}, y \big\rangle - \frac{1}{2s} \|y - y_k \|^2 \right\}.  \label{eqn: ah-ascent}
\end{align}
\end{subequations}
However, the Arrow-Hurwicz algorithm~\eqref{eqn: arrow-hurwicz} suffers from a major drawback: it is not guaranteed to converge generally. This lack of convergence is illustrated in~\Cref{fig: s1}. 
To guarantee the convergence for the problem~\eqref{eqn: composite-gen-lasso},~\texttt{PDHG} is proposed in \cite{chambolle2011first}, which introduces an extrapolated momentum step as: 
\begin{subequations}
\label{eqn: pdhg}
\begin{empheq}[left=\empheqlbrace]{align}
  & x_{k+1} = \argmin_{x \in \mathbb{R}^{d_1}} \left\{ f(x) + \big\langle Fx, y_k \big\rangle + \frac{1}{2s} \|x - x_k \|^2 \right\},                                            \label{eqn: pdhg-descent} \\
  & \overline{x}_{k+1} = x_{k+1} + (x_{k+1} - x_k),                                                                                                                                                                        \label{eqn: pdhg-special} \\
  & y_{k+1} = \argmax_{y \in \mathbb{R}^{d_2}} \left\{ - g^{\star}(y) + \big \langle F\overline{x}_{k+1}, y \big\rangle - \frac{1}{2s} \|y - y_k \|^2 \right\}.     \label{eqn: pdhg-ascent}
\end{empheq}
\end{subequations}

In \cite{li2024understanding}, the high-resolution ODE system for~\texttt{PDHG}~\eqref{eqn: pdhg} is derived as
\begin{subequations}
\label{eqn: pdhg-high}
\begin{align}
  & \dot{X} - s F^{\top} \dot{Y} =-  F^{T}Y - \nabla f(X),                  \label{eqn: high-descent} \\
  & \dot{Y} - s F\dot{X} = FX - \nabla g^{\star}(Y),                          \label{eqn: high-ascent}
\end{align}
\end{subequations}
where the two terms, $s F^{\top} \dot{Y}$ and $ F\dot{X}$, referred to as the combined $x$-correction and $y$-correction, distinguish this system from the one for the Arrow-Hurwicz algorithm~\eqref{eqn: arrow-hurwicz}.  In the Arrow-Hurwicz algorithm,  these two correction terms vanish. As a result, the system degenerates into a Hamiltonian system for the counterexample shown in~\Cref{fig: s1}. This guarantees energy conservation, which prevents convergence. Consequently, the high-resolution ODE system~\eqref{eqn: pdhg-high} successfully explains why the Arrow-Hurwicz algorithm fails to converge in certain scenarios, whereas the \texttt{PDHG} algorithm converges. Similar to the analysis for ADMM, the continuous model enables the construction of a discrete Lyapunov function:
\begin{equation}
\label{eqn: pdhg-lyapunov}
\mathcal{E}(k) = \frac{1}{2s} \|x_k - x^{\star}\|^2 + \frac{1}{2s} \|y_k - y^{\star}\|^2 - \big\langle F(x_k - x^{\star}), y_k - y^{\star} \big\rangle.
\end{equation}
Furthermore, the numerical error induced by implicit discretization can also be interpreted as a secondary Lyapunov function:
\begin{align}
\mathcal{NE}(k) =  & \frac{1}{2s} \|x_{k+1} - x_{k}\|^2 +  \frac{1}{2s}  \|y_{k+1} - y_{k}\|^2 \nonumber \\
                              & - \big\langle F(x_{k+1} - x_{k}), y_{k+1} - y_{k} \big\rangle. \label{eqn: pdhg-numerical-error}
\end{align}
Together, the two Lyapunov functions,~\eqref{eqn: pdhg-lyapunov} and~\eqref{eqn: pdhg-numerical-error}, establish the convergence results for~\texttt{PDHG} analogous to those shown in~\Cref{thm: admm-convegence-rate}. 

% diminish, thus leads to that for the counterexample shown in~\eqref{fig: s1} degenerates to the Hamilton system, which guarantees the energy conservative, thus does not converge.  Thus, based on the high-resolution ODE system~\eqref{eqn: pdhg-high}, the does not convergence phenomenon for the Arrow-Hurwicz algorithm, while the~\texttt{PDHG} algorithm converges is explained. 

An accelerated variant of~\texttt{PDHG} is also proposed in \cite{chambolle2011first}.  By utilizing the iteration-varying parameters $\tau_k, \sigma_k > 0$ satisfying the coupling condition $\tau_{k+1}\sigma_k = s^2$, alongside the extrapolation factor $\theta_k = \tau_{k+1}/\tau_k \in (0, 1)$, the resulting iteration scheme is given by: 
\begin{subequations}
\label{eqn: acc-pdhg}
	\begin{align}
		& x_{k+1} = \arg\min_{x \in \mathbb{R}^{d_1}} \left\{ f(x) + \big\langle Fx, y_k \big\rangle + \frac{1}{2\tau_k} \|x - x_k \|^2 \right\},                                           \label{eqn: pdhg-descent} \\
		& \overline{x}_{k+1} = x_{k+1} + \theta_k (x_{k+1}-x_k),                                                                                                                                                            \label{eqn: pdhg-special} \\
		& y_{k+1} = \arg\max_{y \in \mathbb{R}^{d_2}} \left\{ - g^{*}(y) + \big \langle F\overline{x}_{k+1}, y \big\rangle - \frac{1}{2\sigma_k} \|y - y_k \|^2 \right\},      \label{eqn: pdhg-ascent}
	\end{align}
\end{subequations}
which achieves an accelerated convergence rate of: 
\begin{equation}
\label{eqn: acc-pdhg-rate}
\| x_k - x^{\star} \|^2 \leq \mathrm{O}\left( \frac{1}{k^2} \right).
\end{equation}
In \cite{zeng2025lyapunov}, we design a novel discrete Lyapunov function as
\begin{align}
	\mathcal{E}(k) = \frac{\|x_k - x^{\star}\|^2}{2\tau_k^2} & + \frac{\|y_{k-1} - y^{\star}\|^2}{2s^2}  \nonumber \\	                                                                                 
	             & + \frac{\left\langle F(x_k - x_{k-1}), y_{k-1} - y^{\star} \right\rangle}{\tau_{k-1}} + \frac{\|x_k - x_{k-1}\|^2}{2\tau_{k-1}^2},   \label{eqn: lyapunov-acceleration}
\end{align}
which yields a concise and elegant proof establishing the convergence rate in~\eqref{eqn: acc-pdhg-rate}.

\section{Conclusion and future directions}
\label{sec: conclusion}

In this review, we have provided a comprehensive survey of recent developments in gradient-based optimization, with a particular focus on the application of high-resolution ODE frameworks and newly developed discrete Lyapunov functions. Our analysis has encompassed several key areas: the acceleration mechanism of~\texttt{NAG-SC}; the core dynamics of~\texttt{NAG} and its proximal generalization,~\texttt{FISTA}, including their practical utility in gradient norm acceleration for image deblurring;  the establishment of new results regarding the linear convergence of these methods for strongly convex functions; and the systematic design of novel Lyapunov functions to extend these analytical frameworks to their monotonic counterparts. Beyond these basic  methods, we have demonstrated the utility of high-resolution ODE modeling and Lyapunov analysis in uncovering and explaining the convergence mechanisms of advanced gradient-based algorithms, such as \texttt{ADMM}, \texttt{PDHG}, and its accelerated variants.

%the acceleration mechanism of~\texttt{NAG-SC},~\texttt{NAG} and its proximal generalizations, and the application of~\texttt{FISTA} to image deblurring, alongside a new result for linear convergence for strongly convex functions. Furthermore, we have explored the design of novel Lyapunov functions to extend these frameworks to their monotonic counterparts. Beyond these fundamental methods, we have demonstrated the utility of high-resolution ODE modeling and Lyapunov analysis in uncovering and explaining the convergence mechanisms of advanced algorithms, such as \texttt{ADMM}, \texttt{PDHG}, and their accelerated variants.

% comprehensive review the recent development for the gradient-based optimization method using the high-resolution ODE frameworks, the new developed discrete Lyapunov functions. including the acceleration mechanism of~\texttt{NAG-SC}, \texttt{NAG} and its proximal generalization,~\texttt{FISTA} in image deblurring problem, Linear convergence on strongly convex functions, and design new Lyapunov functions to extend the framework to the monotonic counterparts, moreover, using the high-resolution ODE framework and Lyapunov analysis to understand the mechanisms of advanced gradient-based algorithms, such as~\texttt{ADMM}, \texttt{PDHG} and its accelerated variants. 

While the ODE framework is well-established for minimization problems, its extension to the minimax setting remains incomplete. Such an extension would be highly
valuable for a broad range of applications, including constrained convex optimization, which is often formulated as large-scale saddle-point problems~\cite{benzi2005numerical}, inverse problems and variational data assimilation~\cite{le1986variational}, and classical control challenges such as the linear quadratic regulator (LQR) problem~\cite{anderson2007optimal}. At present, these problems are often solved using the preconditioned generalized minimal residual (GMRES) method~\cite{benzi2004preconditioner}. However, GMRES suffers from increasing memory and computational costs as the number of iterations grows, and its practical performance depends critically on the availability of high-quality preconditioners that effectively exploit the structure of the underlying problem. From an optimization perspective, these challenges are often addressed via the augmented Lagrangian method; yet, its slow convergence necessitates the development of acceleration techniques. Recently, it is demonstrated in~\cite{liu2026accelerated} that incorporating a proximal operator into the augmented Lagrangian can be interpreted as an implicit gradient descent-ascent (\texttt{GDA}) scheme, which naturally facilitates variable step-size extensions. Furthermore, the Lyapunov techniques specifically designed for minimax objective functions in~\cite{liu2026accelerated} offer a promising new avenue for the design and analysis of efficient algorithms for large-scale saddle-point problems, particularly in nonsmooth settings.

\bibliographystyle{plain}
\bibliography{sigproc}

%  The bibliography
%
%\begin{thebibliography}{9}
%%  Use \bibitem{r1} or \bibitem[Surname(2010)]{r1} (for authoryear case)
%
%\bibitem{nesterov1983method}
%Y.~Nesterov.
%\newblock A method of solving a convex programming problem with convergence rate
%  $O(1/k^2)$.
%\newblock {\em Soviet Mathematics-Doklady}, 27(2):372--376, 1983.
%
%\end{thebibliography}

\end{document}